\documentclass[a4paper,11pt]{article}

\usepackage{latexsym}
\usepackage{theorem}
\normalbaselines
\def\ifm#1#2{\relax\ifmmode#1\else#2\fi} 

\newcommand{\xon}    {\ifm {X_1,\ldots,X_n} {$X_1,\ldots,X_n$}}

\newcommand{\fop}    {\ifm {f_1,\ldots,f_p} {$f_1,\ldots,f_p$}}

\newcommand{\Qxon}   {\ifm {\Q[X_1,\ldots,X_n]} {$\Q[X_1,\ldots,X_n]$}}

\newcommand{\klk}    {\ifm {,\ldots,} {$,\ldots,$}}
\newcommand{\plp}    {\ifm {+\cdots+} {$+\ldots+$}}

\newcommand{\be}{\begin{equation}}
\newcommand{\ee}{\end{equation}}

\def \R{{\rm I\kern -2.2pt R\hskip 1pt}} 

\def \N{{\rm I\kern -2.1pt N\hskip 1pt}} 

\def \P{{\rm I\kern -2.2pt P\hskip 1pt}} 
\newcommand{\kpk}{, \ldots ,}

\newcommand{\Strichq}  {\vrule height0.65em width0.05em depth0em \,}
\newcommand{\Q}        {\ifm {\mbox{\rm Q}\hspace{-0.54em}\Strichq\>\;}
                             {$\mbox{\rm Q}\hspace{-0.54em}\Strichq\>\;$}}
\newcommand{\C}        {\ifm {\mbox{\rm C}\hspace{-0.45em}\Strichq\>}
                             {$\mbox{\rm C}\hspace{-0.45em}\Strichq\>$}}
\newenvironment{prf}{\vs\noindent \textbf {Proof}\\}{$\mbox{}$\hfill $\Box$ \\}

\def\cqfd{\vbox{\hrule height 5pt width 5pt }\bigskip} 
\def\qed\cqfd
\def\noi{\noindent}
\def\vs{\smallskip}

\theoremstyle{break}
\newtheorem{proposition}{Proposition}
\newtheorem{definition}[proposition]{Definition} 
\newtheorem{notation}[proposition]{Notation} 
\newtheorem{theorem}[proposition]{Theorem}
 
\newtheorem{lemma}[proposition]{Lemma}
\newtheorem{remark}[proposition]{Remark}

\begin{document}

\title{\bf {Polar Varieties and Efficient Real Elimination}$^{1}$}

\author{\sc B. Bank $^{2}$, M. Giusti $^{3}$, J. Heintz $^{4}$,\\[.5cm] \sc
                  G. M. Mbakop $^{2}$ \\[.5cm] \em
                    Dedicated to Steve Smale}
                 
\maketitle 

\addtocounter{footnote}{1}\footnotetext{Research partially supported by 
the following German, French, Spanish and Argentinian grants:
BA 1257/4-1 (DFG), ARG 018/98 INF (BMBF), UMS 658, ECOS A99E06,
DGICYT PB96-0671-C02-02, ANPCyT 03-00000-01593, UBACYT TW 80 and 
PIP CONICET 4571/96. The first two authors wish to thank the MSRI at Berkeley
for its hospitality during their stay, fall 1998.}
\addtocounter{footnote}{1}\footnotetext{Humboldt-Universit\"at zu Berlin,
Institut f\"ur Mathematik,
10099 Berlin, Germany.
bank@mathematik.hu-berlin.de, mbakop@mathematik.hu-berlin.de}
\addtocounter{footnote}{1}\footnotetext{UMS MEDICIS, Laboratoire GAGE,
\'Ecole Polytechnique, 91228 Palaiseau Cedex, France.
giusti@gage.polytechnique.fr}

\addtocounter{footnote}{1}\footnotetext{Departamento de Matem\'aticas,
Estad\'{\i}stica y Computaci\'on, Facultad de Ciencias, Universidad de
Cantabria, 39071 Santander, Spain. 
heintz@matesco.unican.es  \newline and Departamento de Matem\'atica,
Universidad de Buenos Aires, Ciudad Univ., Pab.I, 1428 Buenos Aires,
Argentina. joos@mate.dm.uba.ar}

\begin{abstract}
\noindent Let $S_0$ be a smooth and compact real variety given by
a reduced regular sequence of polynomials $f_1 \klk f_p$. 
This paper is devoted to the algorithmic problem of finding {\em efficiently}
a representative point for each connected component of $S_0$ .
For this purpose we exhibit explicit 
polynomial equations that describe the generic polar 
varieties of $S_0$. This leads to a procedure which
solves our algorithmic problem in time that is polynomial in the
(extrinsic) description length of the input equations $f_1 \klk f_p$
and in a suitably introduced, intrinsic geometric parameter,
called the  {\em degree} of the real interpretation of the given equation
system  $f_1 \klk f_p$.
\end{abstract}

\smallskip


\section{Introduction}
The core of this paper consists in the exhibition of a system of canonical 
equations which describe locally    
the generic polar varieties of a given {\em semialgebraic} complete
intersection manifold $ S_0$ contained in the real
$n$--dimensional affine space $\R^n$. This purely mathematical
description of the polar varieties allows the design of a new type of efficient 
algorithm 
(with intrinsic complexity bounds), which computes, in case that $S_0$ is
smooth and compact, at least one representative point for each connected
component of  $S_0$ (the algorithm returns each such point in a suitable
symbolic codification).  This new algorithm (and, in particular, its
complexity) is the main practical outcome of the present paper. Let us now
briefly describe  our results.
\smallskip

Suppose that the real
variety $S_0$ is compact and given by
polynomial equations of the following form:
\[
f_1(X_1, \ldots, X_n)= \cdots = f_p(X_1, \ldots , X_n) = 0,
\]
where $p, n \in \N,\; p \le n$ and $f_1 \klk f_p$ belong to the
polynomial ring $\Qxon$ in the indeterminates
$X_1 \klk X_n$ over the rational numbers $\Q$.
Let $d$ be a given natural number and assume that for
$1 \le k \le p$ the total degree $\deg f_k$ of the polynomial
$f_k$ is bounded by $d$.
Moreover, we suppose that the polynomials $f_1 \klk f_p$ form a regular 
sequence in $\Qxon$ and that they
are given by a division-free arithmetic circuit of size $L$ that evaluates
them in any given point of the real (or complex) $n$--dimensional affine
space $\R^n$ (or $\C\;^n$). Further, we assume that the Jacobian 
$J(f_1 \klk f_p)$ of the equation system $f_1= \cdots = f_p=0$ has
maximal rank in any point of $S_0$ (thus, implicitly, we assume that
$S_0$ is smooth). Let $W_0:=V(f_1 \klk f_p)$ denote the
({\em complex}) algebraic variety defined by
the polynomials $f_1 \klk f_p$ in the affine space $\C^n$. 
We denote the singular locus of $W_0$ by $Sing W_0$. 
\smallskip

Moreover, let us suppose that the variables $X_1 \klk X_n$ are in generic
position with respect to the equation system $f_1 \klk f_p$. For 
$1 \le i \le n-p$ let $W_i$ be the $i$--th {\em formal (complex) 
polar variety} associated with $W_0$ (and the variables $X_{p+i} \klk X_n$).
\smallskip

Further, let us denote the real counterpart of 
$W_i$ by $S_i:= W_i \cap \R^n$. We call $S_i$ the $i$--th formal
{\em  real} polar variety associated with the real semialgebraic variety
$S_0$ (and the variables $X_{p+i} \klk X_n$) . It turns out that the
(locally) closed sets $W_i \setminus SingW_0$  (resp. $S_i$) are either empty
or  complex (resp. real) smooth manifolds of dimension $n-(p+i)$. Moreover,
for $1 \le i \le n-p$, one sees easily that
\[\widetilde{W_i}:= \overline{W_i \setminus Sing W_0}\]
is the $i$--th polar variety (in the usual sense) 
associated with $W_0$ and the variables $X_{p+i} \klk X_n$ (here, 
$\overline{W_i \setminus Sing W_0}$ denotes the $\Q\mbox{-Zariski}$ closure in
$\C^n$ of the quasi--affine variety $W_i \setminus Sing W_0$). For a precise
definition of the notion of  formal polar varieties and of polar varieties in
the usual sense we refer to Section 2.
\smallskip

Suppose that the real variety $S_0$ is non--empty and satisfies our
assumptions. In Theorem \ref{th:R} of this paper we show that every real polar variety
$S_i=W_i \cap \R^n,\;1 \le i \le n-p,$ is a non--empty, smooth manifold
of dimension $n-p-i$ containing at least one point of 
each connected component of the real variety $S_0$. In particular,  
the real variety $S_{n-p}$ is a finite set containing at least one
representative point of each connected component of $S_0$.
\smallskip

Under the same assumptions  we show in Theorem \ref{th:T} that for $\;1\le i \le n-p$
the quasi--affine variety $W_i \setminus Sing W_0$ is a locally
complete  intersection that satisfies the Jacobian criterion. More
precisely, the quasi--affine variety $W_i \setminus Sing W_0$ is a smooth
manifold of codimension $p+i$ that can be described locally 
by certain regular sequences consisting of the polynomials $f_1 \klk f_p$ and
$i$ many  well--determined $p$--minors of the Jacobian $J(f_1 \klk f_p)$
of the $f_1 \klk f_p$. 
In particular, the quasi--affine variety $W_{n-p} \setminus Sing W_0$ is 
zero-dimensional, whence $\widetilde{W}_{n-p}=W_{n-p}\setminus Sing W_0$. Thus
$\widetilde{W}_{n-p}$ is a zero-dimensional complex variety that
contains a representative point of each connected component of 
the real variety $S_0$. 
\smallskip

The practical outcome of Theorem \ref{th:T} and  Theorem \ref{th:R} consists in
the design of an efficient algorithm (with intrinsic complexity bounds), which
adapts the elimination procedure for complex algebraic
varieties  developed in \cite{gh2} and \cite{gh3} to the real case.
Under the additional assumption that for any $1 \le k \le p$, the
intermediate ideal $(f_1 \klk f_k)$ generated by $f_1 \klk f_k$ in $\Qxon$ 
is radical, we shall apply this procedure
to the $p{n \choose {p-1}}$ well--determined
equation systems of Theorem \ref{th:T}, which describe the zero-dimensional
algebraic variety $\widetilde{W}_{n-p}= W_{n-p} \setminus Sing W_0$ locally. 
In order to find at 
least one representative
point for every connected component of the real variety $S_0$, 
we have just  to run the procedure of \cite{gh2} and \cite{gh3}
on all these equation systems.
Counting arithmetic operations in $\Q$ at unit costs, this can be done in
sequential time
\[
{n \choose {p-1}}L(nd \delta)^{O(1)},
\]
where $\delta$ is the following geometric invariant of the  regular sequence
$f_1 \klk f_p$:  
\[\delta:=\max\{ \max\{ \deg
\overline{V(f_1 \klk f_k) \setminus Sing W_0} | 1 \le k \le p \},\]
\[\max\{ \deg \widetilde{W_i}|1 \le i \le n-p \}\}\]
(here, $\deg \overline{V(f_1 \klk f_k) \setminus W_0}$ and
$\deg \widetilde{W_i}$ denote the geometric degree in the sense of \cite{he} of
the corresponding algebraic varieties).
\smallskip

This is the content of Theorem \ref{th:C} below.  For any $1 \le k \le p$ and 
any $1 \le i \le n-p$ the quantity $\delta$ bounds the degree of the
algebraic variety $\overline{V(f_1 \klk f_k) \setminus Sing W_0}$ and of the
$i$--th polar variety $\widetilde{W_i}= \overline{W_i \setminus Sing W_0}$.
\smallskip

In \cite{gh2} and \cite{gh3} the quantity 
$\max\{ \deg V(f_1 \klk f_i)| 1 \le i \le p \}$ is called the geometric
degree (of the {\em complex} interpretation) of the equation system $f_1 \klk f_p$. 
In analogy to this terminology,
we shall call $\delta$ the {\em geometric degree} of the {\em real} interpretation
of the equation system $f_1 \klk f_p$.
In view of the complexity result above we shall understand the
parameter  $\delta$ as an {\em intrinsic} measure for the size of the real
interpretation of the given polynomial equation system. 
Nevertheless, the word "intrinsic" should be interpreted with some caution in
this context: observe that the complexity parameter $\delta$ depends rather on
the equations $f_1 \klk f_p$ and their order than just on the variety
$\overline{W_0 \setminus Sing W_0}$. 
\smallskip

In order to make our
complexity result more  transparent we are going now to exhibit, 
in terms of extrinsic parameters, some estimations for the
intrinsic system degree $\delta$.
\smallskip

Let us write $d_1:= \deg f_1 \klk d_p:= \deg f_p$ and let 
$D:= d_1 \cdots d_p$ denote the classical B\'ezout number of the
polynomial system  $f_1 \klk f_p$. Then we have the following  degree
estimations for the complex algebraic variety $W_0=V(f_1 \klk f_p)$
\[
\deg \overline{S_0} \le \deg W_0 \le D \le d^p
\]
($\overline{S_0}$ denotes again the $\Q\!\mbox{-Zariski}$ closure in $\C^n$ of
the  real variety $S_0$ ).
\smallskip

On the other hand, we conclude from Theorem \ref{th:T}  that, for every $i,\; 1 \le i
\le n-p,$ the
polar variety $\widetilde{W_i}$ is defined by the
initial system $f_1 \klk f_p$ and certain $p$--minors of the Jacobian
$J(f_1 \klk f_p)$.
Let us denote  the maximum degree of these
$p$--minors by $c_i$. It turns out
that for, any $1 \le i \le n-p$, the polar variety $\widetilde{W_i}$ is a 
codimension one subvariety of $\widetilde{W}_{i-1}$. Now one sees easily that 
the quantity $D_i:=D\;c_1 \cdots c_i$ represents a reasonable 
"B\'ezout number" of  the variety $\widetilde{W_i}$ and that this 
B\'ezout number satisfies the estimate $\deg \widetilde{W_i} \le D_i$.
Putting all this together, we deduce the following estimate for
the intrinsic system degree $\delta$:
\[
\delta \le D_{n-p}= D c_1 \cdots c_{n-p}.
\]
Observing  that for any $i,\;1 \le i \le n-p,$ the inequality
$c_i \le d_1 + \cdots + d_p -p$ holds,
we find the estimations:
\[
\delta \le D (d_1 + \cdots + d_p -p)^{n-p}
\le d^p(pd -p)^{n-p} <p^{n-p}d^n.
\]
In conclusion, our new real algorithm has a time complexity that is, 
in worst
case, polynomial in the "B\'ezout number" $D c_1\cdots c_{n-p}$ of the
zero--dimensional polar variety $\widetilde{W}_{n-p}$.
\smallskip

Our complexity bound ${n \choose {p-1}}L(nd \delta)^{O(1)}$ depends
 on the intrinsic (geometric, semantic) 
parameter $\delta$
and on the extrinsic (algebraic) parameters $d$ and $n$
in a {\em polynomial} manner, and it depends
 on the syntactic parameter $L$ only {\em linearly}.
In this sense one may consider our complexity bound as {\em intrinsic}.
Our real algorithm promises therefore to be practically applicable to special
equation systems with low value for the intrinsic parameter $\delta$.
\smallskip

On the other hand, {\em even in worst case} our algorithm {\em improves}
upon the known $d^{O(n)}$--time procedures for the  algorithmic problem under
consideration, also in their most efficient versions \cite{basu}, \cite{basu1}
(see also \cite{bar}, \cite{bcss}, \cite{CuSm},
\cite{hroy}, \cite{hroy1}, \cite{hroy2}, \cite{Renegar1}, \cite{rene},
\cite{CaEm95},
\cite{canny}). However, this distinction does not become apparent when we measure
complexities simply in terms of $d$ and $n$ (all  mentioned algorithms have  
worst--case complexities of type $d^{O(n)}$), but it becomes clearly 
visible when we use the "B\'ezout number" just introduced as complexity parameter. 
Only our new algorithm is polynomial in this quantity. On the other hand, we
are only able to reach our goal of algorithmic efficiency by means of a
strict limitation to a purely geometric point of view. For the moment there is
no hope that any of the standard questions of real algebra (e.g. finding
generators for the real radical of a polynomial ideal or the formulation of an
effective real Nullstellensatz) can be solved within the complexity framework
of this paper (compare \cite{lom}, \cite{ben} and \cite{bes}). 
\smallskip

In conclusion we may say that this paper establishes a new connection between
the algorithmic complexity of finding a representative set of real solutions of
a given polynomial equation system and the geometry of the (complex) algebraic
variety defined by this system. However, there is a price to pay for that:
this connection becomes only visible if we restrict ourselves to reduced
complete intersection systems that define smooth, compact real varieties.
\smallskip

Our (algorithmic and mathematical) methods and results represent a non--obvious
generalization of the main outcome of \cite{bank1},
where an intrinsic type algorithm was designed for the problem of 
finding at least one representative point in each connected component of a real,
compact {\em hypersurface} given by an $n$--variate, smooth polynomial equation 
$f$ of degree $d \ge 2$ with rational coefficients (such that $f$
represents a regular equation of that hypersurface). This is the particular case
of codimension $p=1$ of the present paper, and our setting leads to
the complexity bound $L(nd \delta)^{O(1)}$ proved in \cite{bank1}.

\section{Polar Varieties}

\subsection{Notations, Notions and General Assumptions}
Let $X_1\kpk X_n$ be indeterminates (or variables) over the rational numbers
$\Q$ and let polynomials $f_1\kpk f_p \in \Qxon$ with $1 \le p \le n$ be given.
Let $\C^n$ and  $\R^n$ denote the $n$--dimensional affine space over the complex and 
the real numbers, respectively. We think $\C^n$ to be equipped with the
$\Q\!\mbox{-Zariski}$ topology, whereas, on $\R^n$, we consider the strong 
topology.  For any subset $U \subset \C^n$ we denote its
$\Q\!\!\mbox{-Zariski--closure}$ by $\overline{U}$.
By $X := (X_1\kpk X_n)$ 
we denote the vector of variables $X_1\kpk X_n$ and by
$x := (x_1\kpk x_n)$ any point of the affine space $\C^n$ or $\R^n$.
We suppose that the polynomials $f_1 \klk f_p$ form a {\em reduced} regular
sequence in $\Qxon$ (here "reduced" means that for any $1 \le k \le p$ the
ideal $(f_1 \klk f_k)$ is radical). The Jacobian of these 
polynomials is denoted by \[ J(f_1\kpk f_p) :=
\left[ \displaystyle\displaystyle\frac{\partial f_k}{\partial X_j}\right]_
{{1 \le k \le p} \atop {1 \le j \le n}}.
\]
For any point $x \in \C^n$ we write 
\[
J(f_1\kpk f_p)(x)  := \left[ \displaystyle\frac{\partial f_k}
{\partial X_j} (x) \right]_{{1 \le k \le p} \atop {1 \le j \le n}}
\]
for the Jacobian of the polynomials $f_1 \klk f_p$ at $x$.
\smallskip

The common complex zeros of the polynomials
$f_1 \klk f_p$ form an affine, $\Q\!\mbox{-definable}$ 
subvariety of $\C^n$, which we denote by 
\[
W_0 := V(f_1\kpk f_p) := \{ x \in \C^n | f_1(x)= \ldots = f_p(x)=0\}.
\]
A point $x \in W_0=V(f_1 \klk f_p)$ is said to be {\em non--singular} or
{\em smooth} (in $W_0$) if the rank of the Jacobian of $f_1 \kpk f_p$ in
$x$ is $p$.
Otherwise $x$ is called a {\em singular} point of $W_0$.
By $Sing W_0$ we denote the set of all singular points of $W_0$.
Since we suppose that the ideal $(f_1 \klk f_p)$ is radical, our notion
of a smooth point coincides with the usual one for algebraic varieties 
by Jacobi's criterion.

\begin{remark}
If  $x \in W_0$  is smooth,
then the hypersurfaces defined by the polynomials 
$f_1 \kpk f_p$ intersect transversally
at the point$x$.
\end{remark}

\begin{definition}
For every $i,\; 1 \le i \le n-p,$ let $\Delta_i$ denote the set of all common 
complex zeros
of all $p$-minors of the 
Jacobian $J(f_1\kpk f_p)$ corresponding to the columns $\{1 \kpk p+i-1\}$. In
other words,
$\Delta_i$ is the determinantal variety defined by all p-minors of the 
submatrix $J_1^{p+i-1}(f_1 \kpk f_p)$ determined by the columns $\{1\klk p+i-1\}$
of the Jacobian $J(f_1\kpk f_p)$.\smallskip

We introduce the affine variety
\[
W_i :=W_0 \cap \Delta_i
\]
associated with the linear subspace of $\C^n$, namely
\[
X^{p+i-1}:= \{ x \in \C^n | X_{p+i}(x)= \ldots = X_n(x)=0\}
\]
and call $W_i$ the {\em $i$--th formal polar variety} of $W_0$.
\smallskip  

By
\[
\widetilde{W_i}:= \overline{W_i \setminus Sing W_0}
\]
we denote the $i$--th {\em polar variety} (in the usual sense) of 
the variety $W_0$.
\end{definition}

\begin{remark}
\begin{itemize}
\item Our definition of polar and formal polar variety depends rather on the
regular sequence $f_1 \klk f_p$ than on the algebraic variety $W_0$. The ad hoc
term "formal polar variety" is only used in this paper for the purpose of
clarification of our subsequent mathematical arguments.
 \item The index $i$  reflects the expected codimension of the polar variety 
$\widetilde{W_i}$ in $W_0$. With respect to the ambient space $\C^n$,  
the  expected  codimension of $\widetilde{W_i}$ is $p+i$ (see Theorem
\ref{th:T} below for a precise statement).
\item According to our notation, the common zeros of all $p$--minors  of
the Jacobian $J(f_1\kpk f_p)$ form the determinantal variety  $\Delta_{n-p+1}$.
Obviously, we have $Sing\;W_0= W_0 \cap \Delta_{n-p+1}  = W_{n-p+1}$.

\item The formal polar varieties $W_i,\;  1 \le i \le n-p,$  constitue a 
decreasing sequence. In particular, we have
\[W_0 \supset W_1 \supset \cdots \supset W_i \supset 
\cdots \supset W_{n-p} \supset W_{n-p+1}=Sing\;W_0.\]
\end{itemize}
\end{remark}

The concept of polar variety goes back to J.--V. Poncelet. Its development
has a long history: Let us mention among others the contributions of F. Severi,
J. A. Todd, S. Kleiman, R. Piene, D. T. L\^e, B. Teissier, J.--P. Henry,
M. Merle ... (see e.g. \cite{pie} and the references quoted there).

\subsection{Local Description of the Determinantal Varieties}

In this subsection we develop a succinct local description of the 
determinantal varieties $\Delta_i,\; 1 \le i \le n-p$.
The following general Exchange Lemma will be our main tool for this description
(this lemma is used in a similar form in \cite{gihn}).
It describes an exchange relation between certain minors
of a given matrix .
\smallskip

Let $A$ be a given  $(p \times n)$-matrix with entries $a_{ij}$
from an arbitrary commutative ring. 
Let $l$ and $k$ be any natural numbers with $l \le n $ and $k \le 
\min \{p,\;l\}$.
Furthermore, let $I_k := (i_1 \klk i_k)$ 
be an ordered sequence of $k$ different elements from the finite set of
natural numbers $\{ 1 \klk l\}$ and let
$M_A \left( I_k \right) := M_A (i_1 \klk i_k)$
denote the $k$-minor of
the matrix $A$ built up by the first $k$ rows and the columns 
$i_1 \klk i_k$. If it is clear by the context what is the matrix
$A$, we shall just write 
$M(i_1 \klk i_k):=M \left(I_k\right):=M_A \left(I_k \right)$.
\smallskip

\begin{lemma}[Exchange Lemma]
As before let a matrix $A$ and natural numbers $l$ and $k$ be given, as well as
two intersecting index sets $I_k = (i_1 \klk i_k)$ and $I_{k-1} = (j_1 \klk j_{k-1})$.
Then, for suitable numbers
$\varepsilon_j \in\{1,-1\}$ with $j \in I_k \setminus 
I_{k-1}$ we have the following identity: 
\[
(\ast ) \hspace{1cm} M \left( I_{k-1}\right) M 
\left( I_k \right) = \sum_{j \in I_k \setminus 
I_{k-1}} \varepsilon_j \;\;  M \left( I_k  \setminus \{ j\} 
\right) M \left( I_{k-1} \cup \{ j\} \right). 
\]
\end{lemma}

\begin{prf} 
Consider the following $((2k-1)\times (2k-1))$-matrix $L$ with entries from
the given matrix $A$:

\renewcommand{\arraystretch}{1.2}
\arraycolsep2mm
\[
L:= \left[ \begin{array}{c|c}
 & L_1(I_k) \\ O & \vdots \\ & L_{k-1} (I_k)  \\
 \hline\\
L_1(I_{k-1})  & L_1(I_k)\\ \vdots & \vdots\\
L_k(I_{k-1}) &  L_k(I_k) \end{array}
\right].
\]
\smallskip

Here, for any $1 \le j \le k,\; L_j\left( I_k \right)$ 
denotes the row vector of length $k$ that we obtain selecting,
from the $j$--th row of the matrix $A$, the $k$ elements placed in the 
columns $I_k=(i_1 \klk i_k)$. Similarly, $L_j\left( I_{k-1} \right)$
is obtained from the $j$--th row of $A$ selecting the $k-1$ elements placed 
in the columns $I_{k-1}=(j_1 \klk j_{k-1})$.\\
Now it is not difficult to verify the identity $(\ast)$ by 
calculating the determinant $\det L$ of the quadratic matrix $L$ via Laplace 
expansion 
in two different ways. First, by expansion of $\det L$ according to the first 
$k-1$ columns of $L$, we obtain the left--hand side of $(\ast)$, disregarding 
the sign.
Expansion of $\det L$ according to the first $k-1$ rows of $L$ 
leads to the right--hand
side of $(\ast)$. This implies the identity $(\ast)$
for an appropriate choice of the signs $\varepsilon_j$,
with $j \in I_k \setminus I_{k-1}$.
\end{prf}

Let $m \in \Qxon$ denote the $(p-1)$--minor of the Jacobian \linebreak
$J(f_1\kpk f_p)$ given by the first $(p-1)$ rows and columns, i.e., let
\[
m := \det \;
\left[ \displaystyle\displaystyle\frac{\partial f_k}{\partial X_j}\right]_
{{1 \le k \le p-1} \atop {1 \le j \le p-1}}.
\]

We consider the determinantal variety  $\Delta_i$ outside of the
hypersurface
\[
V(m):= \{x \in\C^n\;|\;m(x)=0\}
\]
and denote this localization by $(\Delta_i)_m$ ,
i.e., we set
\[
(\Delta_i)_m :=\Delta_i \setminus V(m).
\] 
From now on , for $1 \le i_1 \le \cdots \le i_p \le n$, let us denote by 
\[
M(i_1 \klk i_p) 
\]
the polynomial in $\Qxon$ defined as       
the $p$--minor of the Jacobian $J(f_1\kpk f_p)$ built up by its $p$
rows and the columns $i_1 \klk i_p$. As before, we denote by
\[
M(i_1 \klk i_p)(x) 
\]
the specialization of $M(i_1 \klk i_p)$ in a  given point 
$x \in \C^n$. 

\begin{proposition}\label{prop:m}
Let $1 \le i \le n-p$ be arbitrarily fixed, and let $m$ be the $(p-1)$--
minor defined above.
Then the determinantal variety $\Delta_i$ is locally (i.e., outside of
the hypersurface $V(m)$)  described by the $i$ polynomials
\[
M(1 \klk p-1,p),\; M(1 \klk p-1,p+1) \klk M(1 \klk p-1, p+i-1).
\] 
In other words, we have  
\[
(\Delta_i)_m := \{x \in \C^n |\; m(x) \not=0, M(1 \klk p-1,s)(x)=0, 
s \in\{p \klk p+i-1\}\},
\]
where $M(1 \klk p-1,s)$ denotes, as above, the $p$--minor of the Jacobian  
$J(f_1\kpk f_p)$ 
built up by the first $p-1$ columns and the $s$--th column.
\end{proposition}
\begin{prf}
It suffices to show that
\[
(\Delta_i)_m \supset \{x \in \C^n |\; m(x) \not=0, M(1 \klk p-1,s)=0, 
s \in\{p \klk p+i-1\}\}
\]
holds.\\ 
Let $ x^{\ast} \in \C^n $ be any point satisfying the conditions $m(x^{\ast})
\not=0$ and \linebreak
$M(1 \klk p-1,s)(x^{\ast})=0$ for every $s \in\{p \klk p+i-1\}$.
We have to verify that  
\[
M(i_1 \klk i_p)(x^{\ast})=0
\]
holds for all ordered $p$--tuples $(i_1 \klk i_p)$ of elements of
$\{1 \klk p+i-1\}$.

Applying the Exchange Lemma to $m=M(1 \klk p-1)$ and 
$M(i_1 \klk i_p)$, we deduce the identity
\[
m(x^{\ast}) M(i_1 \klk i_p)(x^{\ast})=
\]
\[
= \sum_{j \in \{i_1 \klk i_p\} \setminus 
\{1 \klk p-1\}} \varepsilon_j \;\;  M \left( \{i_1 \klk i_p\}\setminus \{ j\} 
\right)(x^{\ast}) M (1 \klk p-1,j )(x^{\ast}) 
\]
for suitable numbers $\varepsilon_j \in \{1,-1\}$ with $j \in
\{i_1 \klk i_p\} \setminus \{1 \klk p-1\}$. By assumption we have
$m(x^{\ast}) \not=0$ and 
$M (1 \klk p-1,j )(x^{\ast})=0\;\mbox{for all}\;j \in \{p \klk p+i-1\}$. This 
implies that $x^{\ast}$ belongs to the set $(\Delta_i)_m$.\\
\end{prf}

\begin{notation}
In the sequel we shall simply write $M_j$ for the 
$p$-- minor $M(1 \klk p-1, j)$ given by 
the first $p-1$ columns of $J(f_1 \klk f_p)$ and the column 
$ j \in \{ p \klk n\}$ . 
\end{notation}
\begin{remark}\label{r:m}
\begin{itemize}
\item Proposition \ref{prop:m} implies that the codimension of $\Delta_i$ outside of
the hypersurface $V(m)$ is at most $i$.
\item Proposition \ref{prop:m} holds also for the determinantal variety 
$\Delta_{n-p+1}$ that defines the singular locus $Sing\;W_0 =W_{n-p+1}$ of
the variety $W_0$. Hence, for any point $x^{\ast} \in \C^n$ satisfying
the condition $m(x^{\ast}) \not=0$ and the $n-p+1$ equations
\[
M_j(x^{\ast})=0,\;j \in \{p\klk n\},
\]
the Jacobian 
$J(f_1 \klk f_p)(x^{\ast})$ becomes singular. 
  
\item Replacing the previously chosen
$(p-1)$-minor $m$ by any other
$(p-1)$-minor  
of the Jacobian $J (f_1 \klk f_p)$, the statement of Proposition \ref{prop:m} 
remains true mutatis mutandis.
\end{itemize}
\end{remark}

\subsection{Local Description of the Formal Polar Varieties}
The aim of this subsection is to show the following fact:
\smallskip

Let the variables $\xon$ be in generic position with respect to the polynomials
$f_1 \klk f_p$, and let $\tilde m$ be any $(p-1)$--minor 
of the Jacobian $J(f_1 \klk f_p)$.
In this subsection we are going to show that any formal polar variety 
$W_i,\;1 \le i \le n-p,$  is
a smooth complete intersection variety outside of the
closed set $ Sing W_0 \cup V(\tilde m)$. Moreover, we shall exhibit a
reduced regular sequence describing this variety outside of 
$ Sing W_0 \cup V(\tilde m)$.\smallskip

As in the previous subsection, let $m \in \Qxon$ denote the $(p-1)$--minor of 
the Jacobian $J(f_1 \klk f_p)$ built up by the first $(p-1)$ rows and 
columns.\smallskip

Let $Y_1 \klk Y_n$ be new variables and let $Y:=(Y_1 \klk Y_n)$. 
For any linear coordinate transformation   
$X= AY$, with $A$ being a regular ${(n \times n)}$--matrix, 
we define the polynomials
\[
G_1(Y) := f_1(AY) \kpk G_p(Y) := f_p(AY).
\]
The Jacobian of $G_1 \kpk G_p$ has the form
\[
J(G_1 \kpk G_p ) := \left[ 
\displaystyle\frac{\partial G_k}{\partial Y_j}\right]_{{1 \le k \le p} 
\atop {1 \le j \le n}} =J(f_1 \kpk f_p)A.
\]

Using a similar notation as before, we denote by
\[
\widetilde{M}(i_1 \klk i_p)
\]
the $p$--minor of the new Jacobian $J(G_1 \klk G_p)$ that corresponds to the
columns $1\le i_1< \cdots < i_p \le n$. 

Moreover, we denote by $\widetilde{M}_j$ the $p$-- minor $\widetilde{M}(1 \klk p-1,j)$ 
determined by the fixed first $p-1$ columns of $J(G_1 \klk G_p)$ and 
the column $j \in \{p \klk n\}$. 
\smallskip

For $p\le r,t \le n$ let $Z_{r,t}$ be a new indeterminate.
Using the following regular $(n-p+1) \times (n-p+1)$--parameter matrix 

\renewcommand{\arraystretch}{1.5}
\[
Z:= \left[ \begin{array}{llclllcl}
1 & 0 && 0 &&&\cdots & 0\\
Z_{p+1,p} & 1 &&&&&&\\
\vdots& \vdots &  \ddots & &\mbox{\Huge O}   && &\vdots \\
Z_{p+i-1,p} & Z_{p+i-1,p+1} &\cdots &  1   &&&&\\
Z_{p+i,p}  &  Z_{p+i,p+1} &\cdots & Z_{p+i, p+i-1}  &  1 &&&\\
\vdots & \vdots &  & \vdots & \vdots &\ddots && 0 \\
Z_{n,p} & Z_{n,p+1} & \cdots & Z_{n,p+i-1} & Z_{n,p+i} & Z_{n,p+i+1}&
\cdots  &1
\end{array} \right],
\]
we construct an $(n \times n)$--coordinate 
transformation matrix $A:=A(Z)$, 
which will enable us to prove the statement at the beginning of this
subsection.\smallskip

For the moment, let us fix an index $1 \le i \le n-p$. 
We consider the formal polar variety $W_i$
outside of the hypersurface $V(m)$. Corresponding
to our choice of $i$,
the matrix $Z$ may be subdivided into submatrices as follows:
\[
Z =\left[ \begin{array}{lc}
Z_1^{(i)} & O_{i,n-p-i+1} \\  Z^{(i)} & Z_2^{(i)}
\end{array}
\right]. 
\]
Here the matrix $Z^{(i)}$ is defined as 

\renewcommand{\arraystretch}{1.6}

\[
Z^{(i)} := \left[ \begin{array}{ccc}
Z_{p+i,p} & \ldots & Z_{p+i,p+i-1} \\
\ldots & \ldots & \ldots \\
Z_{np} & \ldots &Z_{n,p+i-1}
\end{array} \right], 
\]
and $Z^{(i)}_1$ and $Z_2^{(i)}$ denote the quadratic lower triangular matrices
bordering $Z^{(i)}$ in $Z$, and $O_{i,n-p-i+1}$ is the $i \times (n-p-i+1)$
zero matrix. 
Let
\[
A:= A(Z) := \left[ \begin{array}{ccc} I_{p-1} & O_{p-1,i} &
O_{p-1,n-p-i+1}\\ O_{i,p-1} & Z_1^{(i)} & O_{i,n-p-i+1}\\
O_{n-p-i+1,p-1} & Z^{(i)} & Z_2^{(i)} \end{array} \right]. 
\]
Here the submatrices $I_r$ and $O_{r,s}$ are unit or
zero matrices, respectively, of corresponding size, and $Z^{(i)},
Z_1^{(i)},$ and $Z_2^{(i)}$ are the submatrices of the parameter matrix $Z$
introduced before. Thus, $A$ is a regular, parameter dependent 
$(n \times n)$--coordinate transformation matrix.
\smallskip

Like the matrix $Z$, the matrix $A$ contains 
\[
s:=\frac{(n-p) (n-p+1)}{2}
\]
parameters $Z_{r,t}$ which we may specialize into 
any point $z$ of the affine space $\C^s$. For such a point $z \in \C^s$
we denote the corresponding specialized matrices
by $A(z)$, $Z_1^{(i)}(z),\;Z_2^{(i)}(z)$ and $Z^{(i)}(z)$. 
\smallskip

We consider now the coordinate transformation given by $X=AY$ with $A=A(Z)$ and 
calculate the Jacobian $J(G_1 \klk G_p)$ with respect to 
the new polynomials $G_1 \klk G_p$ .
Recall that the coordinate transformation matrix $A$ depends on our 
previous choice of the index $1 \le i \le n-p$. 
\smallskip

According to the structure of the coordinate transformation matrix 
$A=A(Z)$ we subdivide the Jacobian $J(f_1 \klk f_p)$ into three submatrices

\renewcommand{\arraystretch}{1.6}
\[
J(f_1 \klk f_p)=\left[ \begin{array}{ccc} U & V & W \end{array} \right],
\]
with
\[
U:=\left[ \displaystyle\displaystyle\frac{\partial f_k}{\partial X_j}\right]_
{{1 \le k \le p} \atop {1 \le j \le p-1}},\;
V:=\left[ \displaystyle\displaystyle\frac{\partial f_k}{\partial X_j}\right]_
{{1 \le k \le p} \atop {p \le j \le p+i-1}},\;
W:=\left[ \displaystyle\displaystyle\frac{\partial f_k}{\partial X_j}\right]_
{{1 \le k \le p} \atop {p+i \le j \le n}}.
\]
From the identity $J(G_1 \klk G_p) = J(f_1 \klk f_p)\;A$ we deduce that our 
new Jacobian is of the form:
\[
J(G_1 \kpk G_p ) = \left[ 
\displaystyle\frac{\partial G_k}{\partial Y_j}\right]_{{1 \le k \le p} 
\atop {1 \le j \le n}} =
\left[ \begin{array}{ccc} U\; & VZ_1^{(i)}+WZ^{(i)} &\; WZ^{(i)}_2 
\end{array} \right].
\]

We are interested in a local description of the $i$--th formal polar variety 
$W_i=W_0 \cap \Delta_i$ outside
of the hypersurface $V(m)$, where $m$ is the fixed
upper left $(p-1)$-- minor of the Jacobian $J(f_1 \klk f_p)$ (and also of its 
submatrix $U$).
Since the coordinate transformation $X=AY$ leaves the submatrix $U$ unchanged, 
the $(p-1)$--minor $m$ remains fixed under this transformation. 
From Proposition \ref{prop:m} we know that the localized determinantal variety 
$(\Delta_i)_m$ is described by the $i$ equations
\[
M_p=0 \klk M_{p+i-1}=0,
\]
and by the condition $m \not=0$. The $p$-- minors $M_p \klk M_{p+i-1}$ 
defining these equations are built up by 
the submatrix $\left[\;U \;\;V \right]$ of the Jacobian 
$J(f_1 \klk f_p)$. Under the coordinate transformation $A(Z)$
the matrix $\left[\;U \;\;V \right]$ is changed into the submatrix 
\[
\left[ \;U\;\;VZ_1^{(i)}+WZ^{(i)}\right]
\]
of the Jacobian $J(G_1 \klk G_p)$ and 
the $p$--minors $M_p \klk M_{p+i-1}$ are changed into the $p$-- minors
\[
\widetilde{M}_p \klk \widetilde{M}_{p+i-1}
\]
of the matrix $\left[ \;U\;\;VZ_1^{(i)}+WZ^{(i)}\right]$. This implies 
the matrix identity
\[
(**)\qquad \left[\widetilde{M}_p \klk \widetilde{M}_{p+i-1}\right]
=\left[M_p \klk M_{p+i-1}\right]\;Z_1^{(i)}+
\left[M_{p+i} \klk M_n\right]\;Z^{(i)}.
\]

For the previously chosen index
$1 \le i \le n-p$, the coordinate transformation $X = A(Z)Y$ 
induces the following morphism of affine spaces: 
\[
\Phi_i :\; \C^n \times \C^s \to \C^p \times \C^i,
\]
defined by 
\[
(x,z) \longmapsto \Phi_i(x,z) :=  
\left(\;f_1(x) \kpk f_p(x), \widetilde{M}_p(x,z) \kpk 
\widetilde{M}_{p+i-1}(x,z)\;\right).
\]
Consider an arbitrary point $z \in \C^s$. We denote by $\Delta_i^{z}$
the determinantal subvariety of $\C^n$ defined by all $p$--minors
of the matrix \linebreak
$\left[ \;U\;\;VZ_1^{(i)}(z)+WZ^{(i)}(z)\right]$
(which is  a submatrix of the new Jacobian obtained by specializing the
coefficients of the polynomials $G_1 \klk G_p$ into the point $z \in \C^s$).
Writing $W_i^{z}:= W_0 \cap \Delta_i^{z}$, 
one sees immediately that 
the zero fiber $\Phi_i^{-1}(0)$ of the morphism $\Phi_i$ contains the set
\[
(W_i^{z})_m:= W_0 \cap (\Delta_i^{z})_m. 
\]   
In other words, for any arbitrarily chosen point $z \in \C^s$, the zero fiber
$\Phi_i^{-1}(0)$ of the morphism $\Phi_i$ contains the transformed  formal polar 
variety $W_i^{z}$, localized in the hypersurface $V(m)$ and expressed in the old 
coordinates.
\smallskip

We are going now to analyze the rank of the Jacobian of the morphism $\Phi_i$
in an arbitrary point $(x,z) \in \C^n \times \C^s$ with $x \in (W_i^{z})_m$.
Using the subdivision of the parameter matrix $Z$ into
the parts $Z^{(i)},\;Z_1^{(i)}$ and $Z_2^{(i)}$,
the Jacobian $J(\Phi_i)$ of the morphism $\Phi_i$ can be written symbolically 
as
\[
J(\Phi_i)= \left[\displaystyle\frac{\partial \Phi_i}{\partial X}\quad
\displaystyle\frac{\partial \Phi_i}{\partial Z^{(i)}}\quad
\displaystyle\frac{\partial \Phi_i}{\partial Z_1^{(i)}}\quad
\displaystyle\frac{\partial \Phi_i}{\partial Z_2^{(i)}}
\right].
\]

We have \smallskip

\arraycolsep2mm
\[
\left[\displaystyle\frac{\partial \Phi_i}{\partial X}\;
\displaystyle\frac{\partial \Phi_i}{\partial Z^{(i)}}\right]
=\]\[=\left [ \begin{array}{cccc}
J(f_1 \klk f_p) & O_{p,n-p-i+1}  & \cdots & O_{p,n-p-i+1}  \\
 \ast &  \left[ {{\partial \widetilde{M}_p} \over {\partial Z_{p+i,p}}} \klk
{{\partial \widetilde{M}_p} \over {\partial Z_{np}}} \right]
 & \cdots   & O_{1,n-p-i+1} \\
 \vdots &  \vdots & \ddots & \vdots \\
 \ast &  O_{1,n-p-i+1}   & \cdots  & \left[ {{\partial \widetilde{M}_{p+i-1}}
  \over
{\partial Z_{p+i,p+i-1}}}
 \klk {{\partial \widetilde{M}_{p+i-1}} \over {\partial Z_{n,p+i-1}}} \right]
\end{array} \right ], 
\]
where the columns correspond to the partial derivatives of $\Phi_i$ with
respect to the variables 
\[
\xon, Z_{p+i,p} \klk Z_{n,p} \klk Z_{p+i,p+i-1} \klk Z_{n,p+i-1}
\]
(in this order). 
The entries $O_{r,t}$ denote here zero matrices of corresponding size and 
the row matrices
labeled by  "$\ast$" represent the partial derivatives with respect
to the variables $\xon$ of the minors 
$ \widetilde{M}_p \klk \widetilde{M}_{p+i-1}$. 
These row matrices will be irrelevant for our considerations.
\smallskip

Furthermore, the third submatrix 
$\left[\displaystyle\frac{\partial \Phi_i}{\partial Z_1^{(i)}}\right]$
of $J(\Phi_i)$ can be written as

\arraycolsep2mm
\[
\left [ \begin{array}{cccc}
O_{p,i-1} & O_{p,i-2}  & \cdots & 0  \\
\left[ {{\partial \widetilde{M}_p} \over {\partial Z_{p+1,p}}} \klk
{{\partial \widetilde{M}_p} \over {\partial Z_{p+i-1,p}}} \right]  &  O_{1,i-2}
 & \cdots   & 0 \\
 O_{1,i-1} & \left[ {{\partial \widetilde{M}_{p+1}} \over {\partial Z_{p+2,p+1}}} \klk
{{\partial \widetilde{M}_p} \over {\partial Z_{p+i-1,p+1}}} \right] 
 & \cdots & 0\\
 \vdots &  \vdots & \ddots & \vdots \\
 O_{1,i-1} &  O_{1,i-2}   & \cdots  & \left[ {{\partial \widetilde{M}_{p+i-1}}
  \over
{\partial Z_{p+i-1,p+i-2}}} \right]\\
O_{1,i-1} &  O_{1,i-2}   & \cdots  & 0
\end{array} \right ], 
\]
and the last submatrix 
$\left[\displaystyle\frac{\partial \Phi_i}{\partial Z_2^{(i)}}
\right]$
of $J(\Phi_i)$ is a zero matrix since the $p$-- minors 
$\widetilde{M}_p \klk 
\widetilde{M}_{p+i-1}$ are indepedent of the parameters $Z_{r,t}$ occurring in
the submatrix $Z_2^{(i)}$ of the coordinate transformation matrix $A(Z)$.
\smallskip

Therefore, the Jacobian $J(\Phi_i)$ is of full rank $p+i$ wherever the
submatrix 
\[
\widetilde{J}(\Phi_i):=\left[\displaystyle\frac{\partial \Phi_i}{\partial X}\;
\displaystyle\frac{\partial \Phi_i}{\partial Z^{(i)}}\right]
\]
is of full rank $p+i$. On the other hand, considering the $i$ row matrices contained
in $\widetilde{J}(\Phi_i)$ for $p\le j \le p+i-1$      
\[
\left[{{\partial \widetilde{M}_j} 
\over {\partial Z_{p+i,j}}}
 \klk {{\partial \widetilde{M}_j} \over {\partial Z_{n,j}}} \right],
\]
we see that the representation $(**)$
of the transformed $p$--minors $\widetilde{M}_j$ implies the identity
\[
\displaystyle\left[{{\partial \widetilde{M}_j} 
\over {\partial Z_{p+i,j}}}
 \klk {{\partial \widetilde{M}_j} \over {\partial Z_{n,j}}} \right]
= \left[\;M_{p+i} \klk M_n\right].
\]
Thus, we obtain the representation

\arraycolsep2mm
\[
\widetilde{J}(\Phi_i) =\left [ \begin{array}{cccc}
J(f_1 \klk f_p) & O_{p,n-p-i+1}  & \cdots & O_{p,n-p-i+1}  \\
 \ast & \left[M_{p+i} \klk M_n\right] 
 & \cdots   & O_{1,n-p-i+1} \\
 \vdots &  \vdots & \ddots & \vdots \\
 \ast &  O_{1,n-p-i+1}   & \cdots  & \left[M_{p+i} \klk M_n\right]
\end{array} \right ].
\]
Since all entries of the submatrix  $\widetilde{J}(\Phi_i)$ of the Jacobian
$J(\Phi_i)$ belong 
to the polynomial ring $\Q[X_1\klk X_n]$, we see that the rank
of the matrix $J(\Phi_i)$ in a given point $(x,z) \in \C^n \times \C^s$ 
with $x \in {(W_i^z)}_m$ depends only on the choice of $x$. 
According to our localization outside of the hypersurface $V(m)$,  
let us consider an arbitrary smooth point $\tilde{x}$ of $W_0=V(f_1 \klk f_p)$
satisfying the condition $m(\tilde{x}) \not=0$. Suppose that 
the submatrix $\widetilde{J}(\Phi_i)(\tilde{x})$ is not of 
full rank, i.e., that 
\[
rk\; \widetilde{J}(\Phi_i)(\tilde{x}) < p+i 
\]
holds. This latter inequality is valid if and only if all $p$-minors 
$M_{p+i} \klk M_n $ 
of the Jacobian $J(f_1 \klk f_p)$ vanish at $\tilde{x}$. 
Let 
$\tilde{z} \in \C^s$ be any parameter point such that the pair 
$(\tilde{x},\tilde{z})$ belongs to 
the fiber
$\Phi^{-1}_i (0)$ of the morphism $\Phi_i$. Since 
the $p$-minors $\tilde{M}_p \klk \tilde{M}_{p+i-1}$
of the transformed Jacobian $J(G_1 \klk G_p)$ must vanish at 
$(\tilde{x},\tilde{z}),$
we deduce from $(**)$ that 
\[
\left[0 \klk 0\right]
=\left[M_p(\tilde{x}) \klk M_{p+i-1}(\tilde{x})\right]\;Z_1^{(i)}(\tilde{z})
\]
holds (here $Z_1^{(i)}(\tilde{z})$  denotes again the matrix obtained by
specializing the entries of 
$Z_1^{(i)}$ into the corresponding coordinates of the point
$\tilde{z} \in \C^s$).
Because of the lower triangular form of the regular matrix $Z_1^{(i)}$, 
the latter matrix equation holds if and only if the conditions 
\[
M_{p+i-1}(\tilde{x})= \cdots =M_p(\tilde{x})=0. 
\]
are satisfyied.
Therefore, our assumptions on $\tilde{x}$ and $\tilde{z}$ imply
$m(\tilde{x}) \not= 0$ and $M_p(\tilde{x})= \cdots = M_n(\tilde{x})=0.$
However, by Remark \ref{r:m} this means that 
the Jacobian $J(f_1 \klk f_p) (\tilde{x})$ is singular.
Hence, $\tilde{x}$ is not a smooth point of $W_0$, i.e., 
$\tilde{x} \in Sing\;W_0$, which contradicts our assumption on 
$\tilde{x}$.\smallskip

Now, suppose that we are given a point $(\bar{x},{z}) \in \C^n \times \C^s$ 
that belongs to the 
fiber $\Phi_i^{-1}(0)$. Then $\bar{x}$ belongs to $W_0$. Further, suppose 
that
$\bar{x}$ is a smooth point of $W_0$ outside of the 
hypersurface $V(m)$. Let us consider the Zariski--open 
neighbourhood $\tilde{U}$ of $\bar{x}$ consisting of  
all points $x \in \C^n$ with $m(x) \not=0$ and $rk\;J(f_1 \klk f_p) =p$,
i.e., we consider
\[
\widetilde{U}:=\C^n \setminus 
\left( Sing\;W_0 \cup V(m)\right).
\]
We are going to show that the restricted morphism
\[
\Phi_i :\; \tilde{U} \times \C^s \to \C^p \times \C^i
\]
is transversal to the origin $0 \in \C^p \times \C^i$.\smallskip

In order to see this, consider an arbitrary point $(x,z)$ of  
$\tilde{U} \times \C^s$ that satisfies the equation 
$\Phi_i(x,z)=0$.
Thus, $x$ belongs to $\tilde{U}\cap W_0$ and is, therefore, a 
smooth point of $W_0$, which is outside of the hypersurface
$V(m)$. By the preceding considerations on the rank 
of the Jacobian $J(\Phi_i)$ 
it is clear that $J(\Phi_i)$ has the maximal rank $p+i$ at 
$(x,z)$.
This means that $(x,z)$ is a regular point of $\Phi_i$. Since $(x,z)$ was an 
arbitrary point of $\Phi_i^{-1}(0) \cap\;(\tilde{U} \times \C^s)$, 
the 
claimed transversality has been shown.\smallskip

Now, applying the Weak--Transversality--Theorem of Thom--Sard 
(see e.g. \cite{demaz})
to the diagram
\[
\begin{array}{lcc}
\Phi_i^{-1} (0) \cap (\tilde{U} \times \C^s)  & \hookrightarrow &
\C^n \times \C^s \\
& \searrow  & \qquad \downarrow \\
&  & \qquad \;\; \C^s \end{array} 
\]
one concludes that there is a residual dense set 
$\Omega_i$  of parameters $z \in \C^s$
for which transversality holds.
This implies that, for every fixed $z \in \Omega_i$, 
the transformed and localized formal polar variety
\[
W_i^z \setminus \left(Sing\;W_0 \cup V(m)\right)
\]
is either empty or a smooth variety of codimension $p+i$. This variety
can be described locally by the polynomials  
\[
(***) \quad\quad\quad f_1(X) \klk f_p(X),\; \widetilde{M}_p(X,z) \klk  
\widetilde{M}_{p+i-1}(X,z)
\]
that form a regular sequence outside of $Sing W_0 \cup V(m)$.
Up to now, our considerations concerned only the change of coordinates
for an arbitrarily {\em fixed} 
$1 \le i \le n-p$. However,  
$\Omega := \bigcap_{i=1}^{n-p} \Omega_i$ is a dense residual 
parameter set in $\C^s$ from which we can choose a simultaneous 
change of coordinates 
for all $1 \le i \le n-p$. For every choice $z \in \Omega$ and 
$1 \le i \le n-p$
the transformed formal polar variety $W_i^{z}$ is, outside of the closed set
$Sing W_0 \cup V(m)$,
a smooth complete 
intersection variety described by the (local) regular sequence $(***)$. 
One sees now easily that 
the affine space $\R^s$ contains a non--empty residual dense 
set of parameters $z$ such that the conclusions above apply to
the coordinate transformation $X=A(z)Y$.  
Moreover, $z$ can be chosen from $\Q^s$.
\smallskip

Taking into account Proposition \ref{prop:m} and Remark \ref{r:m},
we deduce the following result from our argumentation:
\smallskip

\begin{theorem} \label{th:T}

Let $W_0 =V(f_1 \klk f_p)$ be a reduced complete 
intersection variety given by polynomials $f_1 \klk f_p$ in $\Qxon$
and suppose that
the variables $\xon$ are in generic position with respect to $f_1 \klk f_p$.
Further, let $m$ be the upper left $(p-1)$--minor of the Jacobian 
$J(f_1 \klk f_p)$. 
Then, every formal polar variety $W_i,\;1 \le i \le n-p,$ 
localized with respect to the closed set
$Sing W_0 \cup V(m)$, is either empty or a  
smooth variety of codimension $p+i$ that
can be described by the equations
\[
f_1 \klk f_p, M_p \klk M_{p+i-1},
\]
where $M_j,\;p \le j \le p+i-1,$ is the $p$-- minor of the 
Jacobian $J(f_1 \klk f_p)$ given by the columns $1 \klk p-1,j$.
Then the polynomials \[f_1 \klk f_p, M_p \klk M_{p+i-1}\] form a regular sequence 
outside of $Sing W_0 \cup V(m)$ .
\end{theorem}

\begin{remark}
Taking into account that the  argumentation on the localization 
with respect to the fixed
$(p-1)$-minor $m$ remains valid mutatis mutandis for any other
$(p-1)$-minor $\tilde m$ of the Jacobian $J (f_1 \klk f_p)$, 
Theorem \ref{th:T} can be restated for any fixed $(p-1)$--minor just by
reordering of columns and rows of the Jacobian $J (f_1 \klk f_p)$.
\end{remark}

\subsection{Existence of Real Points in the Polar Varieties}

Let $\fop \in \Qxon$ be a reduced regular sequence and let again 
$W_0 := V(f_1 \klk f_p)$ be the affine variety defined by $\fop$. Consider 
the real variety $S_0:= W_0 \cap \R^n$ and suppose that
\begin{itemize}
\item[(i)]$S_0$ is nonempty and bounded (and hence compact),
\item[(ii)]the Jacobian $J(\fop)(x)$ is of maximal rank in all points $x$ of
$S_0$ (i.e., $S_0$ is a smooth subvariety of $\R^n$ given by the reduced
regular  sequence $\fop$),
\item[(iii)]the variables $\xon$ are in generic position with respect to 
the polynomials $\fop$.
\end{itemize}

Further, let $C$ be any connected component of the compact set $S_0$, and let
\linebreak  $b:=(a_1 \klk a_{p-1},a_p \klk a_{n-1}, a_n) \in C$ be a locally 
maximal point of the last coordinate $X_n$ in the non--empty compact set 
$C \subset S_0$. Without loss of generality we may assume that the upper left
$(p-1)$--minor $m$ of the Jacobian $J(f_1 \klk f_p)$ does not  vanish in $b$ 
(by our assumptions there must be a $(p-1)$--minor of $J(f_1 \klk f_p)$ not 
vanishing at $b$). In any local parametrization of $S_0$ at $b$ 
the variable $X_n$ cannot be an independent variable, since $X_n$ attains a
local maximum in $b$ ($a_n$ is this local maximum). Hence, without loss of 
generality we may assume that the local parametrization of $S_0$
in $b$ has the following form: there exists an open set ${\cal U} \subset \R^{n-p}$ 
containing 
the point $a:=(a_p \klk a_{n-1})$,  and a continuously differentiable function
\[
\varphi: {\cal U}  \to \R^p\;,
\varphi:=(\varphi_1 \klk \varphi_{p-1},\varphi_n)
\]
such that
\[
x_1= \varphi_1(x_p \klk x_{n-1}) \klk x_{p-1}= \varphi_{p-1}(x_p \klk x_{n-1}),
\]
\[ 
x_n= \varphi_n(x_p \klk x_{n-1})
\]
holds for any $x=(x_p \klk x_{n-1}) \in {\cal U}$.  With 
respect to this local parametrization, the polynomials $f_k,\;1 \le k \le p$,
induce real valued functions of the form:
\[
\tilde{f_k}(X_p \klk X_{n-1}):= 
\]
\[
f_k(\varphi_1(X_p \klk X_{n-1}) \klk \varphi_{p-1}(X_p \klk X_{n-1}),\]\[
X_p \klk X_{n-1},\varphi_n(X_p \klk X_{n-1})) .
\]
For every $1 \le k \le p,$ and every $p \le j \le n-1,$ one has the identity
\be 
\frac{\partial \tilde f_k}{\partial X_j} =  \frac{\partial f_k}{\partial X_j}
+\frac{\partial f_k}{\partial X_1}
\frac{\partial \varphi_1}{\partial X_j} \plp 
\frac{\partial f_k}{\partial X_{p-1}}
\frac{\partial \varphi_{p-1}}{\partial X_j}+\frac{\partial f_k}{\partial X_n}
\frac{\partial \varphi_n}{\partial X_j} = 0
\ee
in the open set ${\cal U}$.

Considering the $(p \times p)$--matrix

\renewcommand{\arraystretch}{2.6}
\[
B:= \left[ \begin{array}{ll}
 \displaystyle\frac{\partial f_1}{\partial X_1}  \quad  \ldots &  
 \displaystyle\frac{\partial f_1}{\partial X_{p-1}} 
 \quad 
  \displaystyle\frac{\partial f_1}{\partial X_n} \\ 
  \dotfill & \dotfill \\
 \displaystyle\frac{\partial f_p}{\partial X_1} \quad \ldots &  
 \displaystyle\frac{\partial f_p}{\partial X_{p-1}} 
 \quad
  \displaystyle\frac{\partial f_p}{\partial X_n} \end{array} \right], 
\]
and observing that $B$ is regular in ${\cal U}$, we obtain from (1) that
\be
- \det B(x)\;  \left[ \begin{array}{c}
 \displaystyle\frac{\partial \varphi_1}{\partial X_j}\\  
 \vdots\\
 \displaystyle\frac{\partial \varphi_{p-1}}{\partial X_j}\\  
 \displaystyle\frac{\partial \varphi_n}{\partial X_j}\end{array} \right] =
 (Adj\; B)(x)\;\left[ \begin{array}{c}
 \displaystyle\frac{\partial f_1}{\partial X_j}(x)\\  
 \vdots\\
 \displaystyle\frac{\partial f_{p-1}}{\partial X_j}(x)\\  
 \displaystyle\frac{\partial f_p}{\partial X_j}(x)\end{array} \right]
\ee
holds for any ${x \in \cal U}$ (here $Adj\;B$ denotes the adjoint matrix 
of the matrix $B$). As $b$ is a locally maximal point of $X_n$, 
we have that 
\[
\displaystyle\frac{\partial \varphi_n}{\partial X_j}\;(a) = 0
\] 
holds for every $p\le j \le n-1$.
Thus, equation (2) implies 
\be
B(n,1)\;(b)\;\displaystyle\frac{\partial f_1}{\partial X_j}\;(b) \plp
B(n,p)\;(b)\;\displaystyle\frac{\partial f_p}{\partial X_j}\;(b) = 0
\ee
for every $p \le j \le n-1$ (here for $1 \le k \le p$ we denote  
the entry of the adjoint matrix
$Adj\;B$ at the cross point of the $k$--th column and the last row
by $B(n,k)$). Taking into account
the particular form of the matrix $B$, the equation system (3) means that 

\renewcommand{\arraystretch}{2.6}
\be
\det \left[ \begin{array}{ll}
 \displaystyle\frac{\partial f_1}{\partial X_1}\;(b)  \quad  \ldots &  
 \displaystyle\frac{\partial f_1}{\partial X_{p-1}}\;(b) 
 \quad 
  \displaystyle\frac{\partial f_1}{\partial X_j}\;(b) \\ 
\dotfill & \dotfill \\
 \displaystyle\frac{\partial f_p}{\partial X_1}\;(b) \quad \ldots &  
 \displaystyle\frac{\partial f_p}{\partial X_{p-1}}\;(b) 
 \quad
  \displaystyle\frac{\partial f_p}{\partial X_j}\;(b) \end{array} \right]\; =
 \; 0 
\ee
holds for every $p \le j \le n-1$. Using our notations for the $p$--minors of
the Jacobian $J(f_1 \klk f_p)$, we reinterprete now the equations (4) as
\[
M_p\;(b)= \ldots =M_{n-1}\;(b)=0.
\]
Since  $m(b) \not=0$ holds by assumption, Proposition \ref{prop:m} implies that $b$ 
belongs to the localized determinantal variety $(\Delta_{n-p})_m$.
Therefore, we have $b \in W_0 \cap (\Delta_{n-p})_m$, i.e., the last
formal polar variety 
$W_{n-p}$ contains the point $b$. 
On the other hand, $b$ is a nonsingular point of $W_0$
and belongs therefore to 
$\widetilde{W}_{n-p}=\overline{W_{n-p} \setminus Sing W_0}$. Thus 
$\widetilde{W}_{n-p}$ is a non--empty set of dimension 
zero that contains the real point $b$ of the arbitrarily chosen 
connected component $C$ of the real variety $S_0$. In particular, 
$b \in \widetilde{W}_{n-p} \cap \R^n \subset W_i \cap \R^n=S_i$ holds
for any $1\le i \le n-p$. 
\smallskip

These considerations imply the following result:

\begin{theorem}\label{th:R}

Let  $W_0 := V(f_1 \klk f_p)$ be as in Theorem \ref{th:T}.
If the real variety $S_0:=W_0 \cap \R^n$ is non--empty, bounded and  
smooth, and if  the variables $X_1 \klk X_n$ are in 
generic position with respect to $f_1 \klk f_p$, then every real formal 
polar variety 
$S_i=W_i \cap \R^n,\;1 \le i \le n-p,$ is a non--empty, smooth manifold
of dimension $n-(p+i)$ and contains at least one representative 
point of each connected
component of the real variety $S_0$.
\end{theorem}

\section{Real Equation Solving}

The geometric results of Section 2 allow us to design  a new efficient 
procedure that finds at least one representative point in each 
connected component of a given smooth, compact, real complete intersection 
variety.
 
This procedure will be formulated  in the algorithmic (complexity) model of 
(division-free) arithmetic circuits and networks 
(arithmetic-boolean circuits) over the rational numbers $\Q$.
 
Roughly speaking, a division-free arithmetic circuit $\beta$ over $\Q$ 
is an algorithmic device that supports a step by step evaluation of certain 
(output) polynomials belonging to
$\Qxon$, say $f_1 \klk f_p$. Each step of $\beta$ corresponds either to an input 
from $X_1 \klk X_n$,  to a constant (circuit parameter) from $\Q$ or to an 
arithmetic operation (addition/subtraction or multiplication). We represent
the circuit $\beta$ by a labelled {\em directed acyclic graph (dag)}. 
The size of this dag measures the sequential time requirements of the evaluation  
of the output polynomials $f_1 \klk f_p$ performed by the circuit $\beta$.
 
A (division-free) arithmetic network over $\Q$ is nothing else but an 
arithmetic circuit that additionally contains  decision gates 
comparing rational values or checking their equality, and selector gates
depending on these decision gates.
 
Arithmetic circuits and networks represent non--uniform algorithms, and the
complexity of executing a single arithmetic operation is always counted at unit cost. 
Nevertheless, by means of well known standard procedures our algorithms will  always be
transposable to the uniform {\em random} bit model and they will be practically 
implementable as well.
All this can be done in the spirit of the general asymptotic complexity 
bounds stated in Theorem \ref{th:C} below.
 
Let us also remark that the depth of an arithmetic circuit (or network) 
measures the {\em parallel} time of its evaluation, whereas its size allows an 
alternative interpretation as "number of processors". In this context 
we would like to emphasize the particular importance of counting only 
{\em nonscalar} arithmetic operations (i.e.,only essential multiplications), 
taking $\Q\!\mbox{-linear}$ operations (in particular,
additions/subtractions)  for cost--free.
This leads to the notion of nonscalar size and depth of a given arithmetic 
circuit or network $\beta$. It can be easily seen that the nonscalar size  
determines  essentially the total size of $\beta$ (which takes into account all operations) 
and that the nonscalar depth dominates the logarithms of degree and height
of the intermediate  results of $\beta$.
\smallskip
 
An arithmetic circuit (or network) becomes a sequential algorithm
when  we play a so--called {\em pebble game} on it. By means of pebble games
we are able to introduce a natural 
space measure in our algorithmic model and along with this, a new,
more subtle sequential time measure. If we play a pebble game on a given 
arithmetic circuit, we obtain a so--called {\em straight line program (slp)}. 
In the same way  we obtain a {\em computation tree} from a given arithmetic network.  
For more details on our complexity model we refer to 
\cite{bue}, \cite{vzG}, \cite{vzG93}, \cite{krick-pardo1}, 
\cite{mat}, \cite{hmw} and especially to
\cite{gls} (where also the implementation aspect is treated).
\smallskip

In the next Theorem \ref{th:C} we are going to consider families of polynomials
$f_1 \klk f_p$ belonging to $\Qxon$, for which we arrange the following 
assumptions and notations:
\begin{itemize}
\item[(i)] $f_1 \klk f_p$ form a regular sequence in $\Qxon$,
\item[(ii)]for every $1 \le k \le p$ the ideal $(f_1 \klk f_k)$ 
generated by $f_1 \klk f_k$ in $\Qxon$ is radical and defines a subvariety of
$\C^n$ of dimension $n-k$ that we denote by $V_k:=V(f_1 \klk f_k)$.
\item[(iii)]the variables $X_1 \klk X_n$ are in
generic position  with respect to the polynomials $f_1 \klk f_p$.
\end{itemize}
Let  
$W_0:= \{x \in \C^n | f_1(x)= \cdots =f_p(x)=0\}$ 
and denote by $Sing W_0$ the singular locus of $W_0$. 
For $1 \le i \le n-p$ let $W_i$ be the $i$--th formal polar variety associated 
with $W_0$ and the variables $X_{p+i} \klk X_n$, and let 
$\widetilde{W_i}:=\overline{W_i \setminus Sing W_0}$ be the $i$--th
polar variety of $W_0$ in the usual sense (see Section 2 for precise 
definitions). Further, for $1 \le k \le p$ 
we shall write  $\widetilde{V}_k:=\overline{V_k \setminus Sing W_0}$.
We call
\[
\delta:=\max\{ \max\{ \deg \widetilde{V}_k| 1 \le k \le p \},
\max\{ \deg \widetilde{W_i}|1 \le i \le n-p \}\}
\] 
the {\em degree} (of the real interpretation) {\em of the polynomial 
equation system}
$f_1 \klk f_p$. Finally, let us make the following assumption:
\begin{itemize}
\item[(iv)]the specialized Jacobian $J(f_1 \klk f_p)(x)$ has maximal rank in 
any point $x$ of
$S_0:=W_0 \cap \R^n=\{x \in \R^n|f_1(x)= \cdots = f_p(x)=0\}$
and $S_0$ is a bounded semialgebraic set
(hence, $S_0$ is empty or a smooth, compact real manifold 
of dimension $n-p$; see Section 2 for details).
\end{itemize}

\begin{theorem}\label{th:C}
Let $n,p,d,\delta,L$ and $\ell$ be natural numbers with $d \ge 2$ and $p \le n$.
There exists an arithmetic network $\mathcal N$ over $\Q$ of size 
${n \choose {p-1}}L(nd\delta)^{O(1)}$ and nonscalar depth 
$O(n(\log nd + \ell)\log \delta)$ with the following property: 
Let $f_1 \klk f_p$ be a family of $n$--variate polynomials
of a degree at most $d$ and assume that $f_1 \klk f_p$ are given by a 
division--free arithmetic circuit $\beta$ in $\Qxon$ of size $L$ and nonscalar 
depth $\ell$. Suppose that the polynomials $f_1 \klk f_p$ satisfy the 
conditions (i), (ii), (iii) and (iv) above. 
Further, suppose that the degree of the real 
interpretation of the polynomial system $f_1 \klk f_p$ is bounded by $\delta$
(let us now freely use the notations just introduced before).
\smallskip

The algorithm 
represented by the arithmetic network $\mathcal N$
starts from 
the circuit $\beta$ as input and decides first whether the complex variety 
$\widetilde{W}_{n-p}$ is empty. If this is not the case, then 
$\widetilde{W}_{n-p}$ is a zero--dimensional complex variety and the
network $\mathcal N$ produces an arithmetic circuit in $\Q$ of 
asymptotically the same size and nonscalar depth as  $\mathcal N$, which 
represents the coefficients of $n+1$ univariate polynomials 
$q,p_1 \klk p_n \in \Q[X_n]$ satisfying the following conditions:
\[ \deg q = \#\;\widetilde{W}_{n-p}, \]
\[\max\{\deg p_k|1\le k \le n\} < \deg q, \]
\[\widetilde{W}_{n-p} = \{(p_1(u) \klk p_n(u))|u \in \C, q(u)=0\}.\]
Moreover, the algorithm represented by the arithmetic network 
$\mathcal N$ decides whether the set $\widetilde{W}_{n-p} \cap \R^n$ is empty. 
In this case we conclude $S_0=W_0 \cap\R^n= \emptyset$. Otherwise, the network
$\mathcal N$ produces at most  $\#\;\widetilde{W}_{n-p} \le \delta$ sign
sequences belonging to the set $\{-1,0,1\}$ such that these sign sequences 
encode the real zeros of the polynomial $q$ "\`a la Thom" (\cite{CosteRoy}).
In this way, namely by means of the Thom encoding of the real zeros of $q$ and
by means of the polynomials $p_1 \klk p_n$, the arithmetic network $\mathcal N$
describes the finite, non--empty set
\[\widetilde{W}_{n-p}\cap \R^n=\{(p_1(u) \klk p_n(u))| u \in \R, q(u)=0\},\]
which contains at least one  representative point for each connected component
of the real variety 
\[S_0=\{x \in \R^n|f_1(x)= \cdots = f_p(x)=0\}.\]
\end{theorem}
\begin{prf}
We shall freely use the notations of Section 2. Any selection of indices 
$1 \le i_1 < \cdots < i_p \le n$ and $1 \le j,k \le p$ determines a 
$p$-- minor $M(i_1 \klk i_p)$ and a $(p-1)$-- minor
$m(i_1 \klk i_p;j,k)$ of the Jacobian $J(f_1 \klk f_p)$ in the
following way: $M(i_1 \klk i_p)$
is the determinant of the $(p \times p)$-- submatrix of $J(f_1 \klk f_p)$
with columns $i_1 \klk i_p$, and $m(i_1 \klk i_p;j,k)$ is the determinant 
of the matrix obtained from the former one deleting the row number $j$ and 
the column number $i_k$. There are 
$p^2{n \choose {p}}$ such possible selections. Let us fix one of them, say 
$i_1:=1 \klk i_p:=p;j:=p,k:=p$. Then, using the notations of Section 2,
we have 
$m(i_1 \klk i_p;j,k)=m,\; M(i_1 \klk i_p)=M_p$. Let us abbreviate $g:=mM_p$.
From our assumptions on $f_1 \klk f_p$ and Theorem \ref{th:T} and 
Theorem \ref{th:R} of Section 2 we 
deduce the following facts: For any $1 \le i \le n-p$ the polynomials
$f_1 \klk f_p, M_p \klk M_{p+i-1}$ have degree at most $pd$. They generate 
the trivial ideal or form a regular sequence in the localized $\Q$-algebra
$\Qxon_g$. In either case the ideal generated by 
$f_1 \klk f_p, M_p \klk M_{p+i-1}$ in
$\Qxon_g$ is radical and defines a complex variety that is empty or of degree
\[\deg (\overline{W_i \setminus V(g)}) 
\le \deg (\overline{W_i \setminus Sing W_0}) =\deg \widetilde{W_i} \le \delta.\]
Moreover, by assumption, the polynomials $f_1 \klk f_p$ form  a regular 
sequence in $\Qxon_g$ and for each $1 \le k \le p$ the ideal generated by 
$f_1 \klk f_k$ in $\Qxon_g$ is radical and defines a complex variety of degree
\[\deg \overline{(V_k \setminus V(g))} \le
\deg \widetilde{V}_k \le \delta.\]
One sees easily that the polynomials $f_1 \klk f_p, M_p \klk M_{n-1}$
and $g$ can be evaluated by a division--free arithmetic circuit of size 
$O(L+n^5)$ and nonscalar depth $O(\log n +\ell)$. Applying now, for each 
$1 \le i \le n-p$, the 
algorithm underlying \cite{gh2}, Proposition 18 in its rational version
\cite{gh3}, Theorem 19 to the system
\[f_1=0 \klk f_p=0,\;M_p=0 \klk M_{p+i-1}=0, g \neq 0\]
we are able to check whether the particular system
\[f_1=0 \klk f_p=0,\;M_p=0 \klk M_{n-1}=0, g \neq 0\]
has a solution in $\C^n$. If this is the case, then this system defines a 
zero--dimensional algebraic set, namely $W_{n-p} \setminus V(g)$, and the algorithm 
produces an arithmetic circuit $\bar \gamma$ in $\Q$ that represents the 
coefficients of $n+1$ univariate polynomials $\bar q,\bar p_1 \klk \bar p_n
\in \Q[X_n]$ satisfying the following conditions:
\[ \deg \bar q = \#\;(W_{n-p} \setminus V(g)),\]
\[\max\{\deg \bar p_k|1\le k \le n\} < \deg \bar q, \]
\[ W_{n-p} \setminus V(g) = \{(\bar p_1(u) \klk \bar p_n(u))|u \in \C,
\bar q(u)=0\}.\] The algorithm is represented by an arithmetic network of size 
$L(nd\delta)^{O(1)}$ and nonscalar depth $O(n(\log nd +\ell)\log \delta)$,
and the circuit $\bar \gamma$ has asymptotically the same size and nonscalar
depth. Running this procedure for each selection $1 \le i_1< \cdots <i_p \le
n$ and $1 \le j,k \le p$ we obtain an arithmetic network $\mathcal N_0$ of size
$p^2{n \choose {p}}L(nd\delta)^{O(1)}={n \choose {p-1}}L(nd\delta)^{O(1)}$ 
and nonscalar depth 
$O(n(\log nd +\ell)\log \delta)$, which decides whether 
$\widetilde{W}_{n-p}=W_{n-p} \setminus Sing W_0$ is empty. Suppose 
that this is not the case. 
Then $\mathcal N_0$ describes locally the variety $\widetilde{W}_{n-p}$, 
which is now zero-dimensional. Each local description of 
$\widetilde{W}_{n-p}$  contains an arithmetic circuit 
representation of the coefficients of the minimal polynomial of  
the variable $X_n$ with respect to the corresponding 
local piece of $\widetilde{W}_{n-p}$ . 
Moreover, one easily obtains the same type of information  
for any linear form $X_i+X_n$  and any variable $X_i$ with $1\le i < n$.
One multiplies now all minimal polynomials of the variable 
$X_n$ obtained in this way. 
Making this product squarefree 
(see e.g \cite{krick-pardo1}, Lemma 12) one obtains the polynomial 
$q$ of the statement of Theorem \ref{th:C}. Doing the same thing for the 
minimal polynomials of each linear form $X_i+X_n$ and each variable 
$X_i$  with $1\le i < n$, yields by means of \cite{krick-pardo1}, 
Lemma 26, the polynomials $p_1 \klk p_n$ of the statement of Theorem \ref{th:C}. 
All this can be done by means of an arithmetic network $\mathcal N_1$, which extends 
$\mathcal N_0$ and has asympotically the same size and nonscalar depth. 
The desired arithmetic network $\mathcal N$ is now obtained
from $\mathcal N_1$ in the  same way as in the proof \cite{bank1},~Theorem~8,
namely as follows: applying the main algorithm of 
\cite{BenOrKozenReif} or \cite{RoySzpirglas} and adding suitable comparison
gates for rational numbers, we extend  $\mathcal N_1$ to a new arithmetical
network $\mathcal N$ of asymptotically the same size and depth, such that 
$\mathcal N$ decides whether the univariate polynomial $q$ has a real zero.
If this is the case, the network $\mathcal N$ enumerates the existing real
zeros of $q$, encoding them "\`a la Thom" (\cite{CosteRoy}). If $q$ has no
real zero we conclude $S_0= \emptyset$. Otherwise, the network $\mathcal N$
encodes all real zeros of $q$ by means of $\# \widetilde{W}_{n-p} \le \delta$
sign sequences belonging to the set $\{-1,0,1\}$. This encoding and the
polynomials $p_1 \klk p_n$ describe now the set $\widetilde{W}_{n-p} \cap \R^n=
\{(p_1(u) \klk p_n(u))|u \in \R,q(u)=0\}$ that contains a representative point
for each connected component of $S_0$.
\end{prf}

\begin{remark}
\begin{enumerate}
\item[(i)]Using the refined algorithmic techniques 
of \cite{hmw} or \cite{gls} it is not too difficult to see 
that for inputs $f_1 \klk f_p$ represented by straight--line programs of 
length $T$ and space $S$ the arithmetic network $\mathcal N$ can be converted into 
an algebraic computation tree which solves the algorithmic problem of 
Theorem \ref{th:C} in time $O((Tdn^2+n^5)\delta^3 \log^3\delta \log^2\log \delta)$ 
and space $O(Sdn\delta^2)$.
\item[(ii)]The smooth,
compact hypersurface case (with $p:=1$) of Theorem \ref{th:C} corresponds exactly to 
\cite{bank1}, Theorem 8.
\item[(iii)]Let $J(f_1 \klk f_p)^T$ denote the transposed matrix of the Jacobian 
\linebreak
$J(f_1 \klk f_p)$ of the polynomials $f_1 \klk f_p$ in the statement of 
Theorem \ref{th:C} and let 
\[\mathcal D:=\det J(f_1 \klk f_p)J(f_1 \klk f_p)^T.\]
From the well--known Cauchy--Binet formula one deduces easily that, with the notations
of Section 2, the identity 
\[\mathcal D= \sum_{1\le i_1< \cdots < i_p \le n}{\det}^2\; M(i_1 \klk i_p)\]
holds. Replacing now, in the statement and the proof of Theorem \ref{th:C} for
$1 \le i \le n-p$, the polar variety $\widetilde{W}_i$ by 
$\widehat{W}_i:=\overline{W_i \setminus V(\mathcal D)}$ and the parameter $\delta$
by
\[
\widehat{\delta}:=\max\{ \max\{ \deg \widetilde{V}_k| 1 \le k \le p \},
\max\{ \deg \widehat{W_i}|1 \le i \le n-p \}\}
\]
one obtains a somewhat improved complexity result, since 
$\widehat{\delta} \le \delta$ holds.
\end{enumerate}
\end{remark}

Let us now suppose that the polynomials $f_1 \klk f_p \in \Qxon$
satisfy the conditions (i), (ii), (iii), (iv) above. Unfortunately, the
complexity parameter $\delta$ of Theorem \ref{th:C} is strongly related to the 
{\em complex} degrees of the polar varieties 
$\widetilde{W}_1 \klk \widetilde{W}_{n-p}$ of 
$W_0=\{x \in \C^n|f_1(x)= \cdots= f_p(x)=0\}$ and not to their {\em real}
degrees. Under
some additional algorithmic assumptions, which we are going to explain below, 
we may replace the complexity parameter $\delta$ by a smaller one that measures 
only the real degrees of the polar varieties 
$\widetilde{W}_1 \klk \widetilde{W}_{n-p}$. We shall call this new complexity 
parameter the {\em real degree} of the equation system  $f_1 \klk f_p$ and denote 
it by $\delta^*$.
\smallskip

Let $1 \le k \le p$ and let us consider the decomposition of the intermediate variety
$\widetilde{V}_k$ into irreducible components with respect to the
$\Q\mbox{-Zariski}$ topology of $\C^n$ say $\widetilde{V}_k=C_1\cup \cdots
\cup C_s$. We call an irreducible  component $C_r,\;1 \le r \le s$, {\em real}
if $C_r \cap \R^n$ contains a smooth  point of $C_r$. The union of all real
irreducible components of $\widetilde{V}_k$ is called the {\em real part} of
$\widetilde{V}_k$ and denoted by $V_k^*$. We call 
$\deg V_k^*$ the {\em real degree} of the intermediate variety
$\widetilde{V}_k$. Similarly, we introduce for every $1 \le i \le n-p$
the real part $W_i^*$ of the polar variety $\widetilde{W}_i$ and its real
degree $\deg W_i^*$. Finally, we
define the {\em real degree of the equation system} $f_1 \klk f_p$ as 
\[\delta^*:=\max\{ \max\{ \deg V_k^*
| 1 \le k \le p \},\max\{ \deg W_i^*|1 \le i \le n-p \}\}.\] 
\smallskip

Now, we are going to restate the main outcome of Theorem \ref{th:C} 
in terms of the new complexity parameter $\delta^*$. For this purpose we have 
to include the following two subroutines in our algorithmic model:

\begin{itemize}
\item the first subroutine we need is a factorization algorithm for
univariate polynomials over $\Q$. In the bit complexity model the problem
of factorizing univariate polynomials over $\Q$ is known to be of
polynomial time complexity \cite{LLL}, whereas in the arithmetic model we
are considering here this question is more intricate \cite{vzGaSe91}.
In the extended complexity model we are going to consider here, the cost of
factorizing a univariate polynomial of degree $D$ over $\Q$ (given by its
coefficients) is accounted as $D^{O(1)}$.
\item  the second subroutine allows us to discard non-real irreducible
components of the occurring complex polar varieties. This second subroutine
starts from a straight-line program for a single polynomial in $\Qxon$
as input and decides whether this polynomial has a real zero
(however, without actually exhibiting it if there is one). Again this
subroutine is taken into account at polynomial cost.
\end{itemize}
We call an arithmetic network  over $\Q$ {\em extended} if it
contains extra nodes corresponding to the first and second subroutine.
\smallskip

Modifying our algorithmic model in this way, we are able to formulate the following 
complexity result, which generalizes \cite{bank1}, Theorem 12 and improves 
the complexity outcome of our previous Theorem \ref{th:C}.
\begin{remark}
Let $n,p,d,\delta^*,L$ and $\ell$ be natural numbers with $d \ge 2$ and $p \le n$.
There exists an extended arithmetic network $\mathcal N^*$ over $\Q$ of size 
${n \choose {p-1}}L(nd\delta^*)^{O(1)}$ with the following property: 
Let $f_1 \klk f_p$ be a family of $n$--variate polynomials
of a degree at most $d$ and assume that $f_1 \klk f_p$ are given by a 
division--free arithmetic circuit $\beta$ in $\Qxon$ of size $L$. 
Suppose that the polynomials $f_1 \klk f_p$ satisfy the 
conditions (i), (ii), (iii), and (iv) contained in the formulation of
Theorem \ref{th:C}. 
Let us now freely use the notations introduced in the present
section. Assume that the real variety
$S_0=\{x \in \R^n|f_1(x)= \cdots = f_p(x)=0\}$ is not empty and that the real
degree of the polynomial system $f_1 \klk f_p$ is bounded by $\delta^*$.
The algorithm represented by the arithmetic network $\mathcal N^*$
starts from the circuit $\beta$ as input and decides first whether the complex variety 
$W_{n-p}^*$ is empty. If this is not the case, then 
$W_{n-p}^*$ is a zero--dimensional complex variety 
and the network $\mathcal N^*$ produces an arithmetic circuit in $\Q$ of 
asymptotically the same size as  $\mathcal N^*$,which 
represents the coefficients of $n+1$ univariate polynomials  
$q^*,p_1^* \klk p_n^* \in \Q[X_n]$ satisfying the conditions 
\[ \deg q^* = \#\;W_{n-p}^*, \]
\[\max\{\deg p_k^*|1\le k \le n\} < \deg q^*, \]
\[W_{n-p}^* = \{(p_1^*(u) \klk p_n^*(u))|u \in \C, q^*(u)=0\}.\]
Each over $\Q$ irreducible component of the complex variety $W_{n-p}^*$
contains at least one real point characterized by an 
irreducible factor of the polynomial $q^*$. The algorithm represented by the network 
$\mathcal N^*$ returns all these points in a codification 
"\`a la Thom". Moreover, the non--empty set $W_{n-p}^* \cap \R^n$
contains at least one 
representative point for each connected component of the real variety
$S_0$.
\end{remark}

The proof of this remark is a straight--forward adaptation of the 
arguments of the proof of \cite{bank1}, Theorem 12 (which treats only the 
hypersurface case with $p:=1$) to the arguments of Theorem \ref{th:C}
above. Therefore, we omit this proof.
\smallskip 

Let us finally observe that the practical relevance of the complexity outcome
of  Remark 6 is highly hypothetical, because it depends on the strong
assumption that extended arithmetical networks are realizable by performant,
programmable algorithms. Nevertheless, by means of Remark 6, we wish to
underline the importance of the search for efficient procedures that realize
the first and second subroutine introduced as extra nodes in our complexity
model of extended arithmetic networks.


\end{document}